\documentclass[11pt,a4paper]{article}
\usepackage[T1]{fontenc}
\usepackage[utf8]{inputenc}
\usepackage{lmodern}
\usepackage[margin=27mm]{geometry}
\usepackage{amsmath,amssymb,amsthm,mathtools}
\usepackage{enumitem}
\usepackage{indentfirst}
\usepackage[colorlinks=true,linkcolor=blue,citecolor=blue,urlcolor=blue]{hyperref}
\allowdisplaybreaks[2]

\newcommand{\NotesAI}{ChatGPT 6.2 Astra (Openai) and Claude Opus 5.5 (Anthropic)}

\newtheorem{theorem}{Theorem}[section]
\newtheorem{proposition}[theorem]{Proposition}
\newtheorem{lemma}[theorem]{Lemma}
\newtheorem{corollary}[theorem]{Corollary}
\newtheorem{mainthm}{Theorem}

\theoremstyle{definition}

\newtheorem{construction}[theorem]{Construction}
\newtheorem{example}[theorem]{Example}
\newtheorem{question}[theorem]{Question}
\theoremstyle{remark}
\newtheorem{remark}[theorem]{Remark}

\DeclareMathOperator{\TC}{TC}
\DeclareMathOperator{\Lie}{Lie}
\DeclareMathOperator{\Fr}{Fr}
\DeclareMathOperator{\Spec}{Spec}
\DeclareMathOperator{\Spf}{Spf}

\DeclareMathOperator{\Tor}{Tor}
\DeclareMathOperator{\rk}{rk}
\DeclareMathOperator{\Frac}{Frac}
\newcommand{\length}{\operatorname{length}}
\newcommand{\Ga}{\widehat{\mathbf G}_a}
\newcommand{\Gm}{\widehat{\mathbf G}_m}
\newcommand{\Wh}{\widehat W}
\newcommand{\Fp}{\mathbf F_p}
\newcommand{\Zp}{\mathbf Z_p}
\newcommand{\Zpl}{\mathbf Z_{(p)}}
\newcommand{\E}{\mathbb E}
\newcommand{\cO}{\mathcal O}
\newcommand{\fm}{\mathfrak m}

\title{Cartier base change and specialization\\ of smooth formal groups}
\author{Yuan Yang\\[2pt]
  \small BICMR, Peking University, Beijing, China\\
  \small\href{mailto:yy\_@bicmr.pku.edu.cn}{\texttt{yy\_@bicmr.pku.edu.cn}}}
\date{\today}

\begin{document}
\maketitle
\vspace{-2em}

\begin{abstract}
Let $G$ be a smooth commutative formal group over a local domain of
characteristic $p>0$. We prove that under specialization the dimensions of
the $p^m$-torsion subgroups of $G$ can only increase, and that the Newton
polygon of the $p$-divisible part, continued by a ray of slope one, can only
rise. The proof compares the orders of the finite group schemes
$\ker(\Fr^n)\cap G[p^m]$ on the two fibres. We then study Cartier modules
over perfect rings. The quotients $M/V^nM$ of the Cartier module commute with
extension of Witt scalars. This gives semicontinuity of the lengths of
$M/(V^nM+p^mM)$ without error terms, and a second proof of the inequalities
for families over $k[[t]]$. The full Cartier module does not commute with
this base change, already for $\Ga$. Finally, we construct explicit families
realizing the specializations $\Gm\times\Wh_r\rightsquigarrow\Wh_{r+1}$ and
$\Wh_a\times\Wh_b\rightsquigarrow\Wh_{a-1}\times\Wh_{b+1}$ for $a\ge b+2$.
For $\Gm\times\Ga\rightsquigarrow\Wh_2$ we give a polynomial formal group law
over $\Fp[t]$.
\end{abstract}

\section{Introduction}

Throughout, $p$ is a prime number. Let $\kappa$ be an algebraically closed
field of characteristic $p$, and let $G$ be a smooth commutative formal group
of finite dimension over $\kappa$. By the structure theory of such groups
\cite[5.36, 5.38]{Zin84} and the Dieudonn\'e--Manin classification
\cite{Man63}, \cite[Ch.~IV]{Dem72}, $G$ is isogenous to a product
\begin{equation}\label{eq:intro-decomposition}
  \Wh_{r_1}\times\cdots\times\Wh_{r_u}\times L_1\times\cdots\times L_v ,
\end{equation}
where $\Wh_r$ is the formal completion of the group of Witt vectors of
length~$r$ and $L_1,\dots,L_v$ are isoclinic $p$-divisible formal groups. We
use the covariant slope convention: an isoclinic $p$-divisible formal group
of height $h$ and dimension $d$ has slope $1-d/h\in[0,1)$, so that $\Gm$ has
slope~$0$ (see \S\ref{subsec:invariants}). The isogeny class of $G$ is
detected by two invariants.
\begin{enumerate}[label=(\roman*)]
\item The dimensions $\dim G[p^m]$, $m\ge1$, of the formal $p^m$-torsion
  subgroups. One has $\dim G[p^m]=\sum_i\min(m,r_i)$, so these dimensions
  determine the multiset $\{r_1,\dots,r_u\}$.
\item The \emph{extended Newton polygon} $P_G\colon[0,\infty)\to\mathbf R$:
  the Newton polygon of $L_1\times\cdots\times L_v$, continued by a ray of
  slope~$1$. Since all slopes of the $L_j$ are smaller than $1$, $P_G$
  determines the Newton polygon of $L_1\times\cdots\times L_v$.
\end{enumerate}
Both invariants depend only on the isogeny class of $G$
(Corollary~\ref{cor:invariants}). Our first result says that both can only
go up under specialization.

\begin{mainthm}\label{thm:A}
Let $\Lambda$ be a local domain of characteristic $p$, and let $G$ be a smooth
commutative formal group over $\Lambda$, given by a formal group law in
finitely many variables. Let $\bar\eta$ and $\bar s$ be geometric points over
the generic point and the closed point of $\Spec\Lambda$. Then for all
$m\ge1$ and all $y\ge0$,
\[
  \dim G_{\bar\eta}[p^m]\le\dim G_{\bar s}[p^m]
  \qquad\text{and}\qquad
  P_{G_{\bar\eta}}(y)\le P_{G_{\bar s}}(y).
\]
\end{mainthm}

For $p$-divisible groups the polygon inequality is Grothendieck's
specialization theorem \cite{Gro74}; see also \cite{Kat79} and
\cite[Ch.~IV, \S7]{Dem72}. For smooth formal groups, the $p$-divisible part
need not form a family. In the family of Theorem~\ref{thm:C}(1) with $r=1$,
the generic fibre is $\Gm\times\Ga$, whose $p$-divisible part $\Gm$ has
height one, while the special fibre $\Wh_2$ is unipotent. The ray of slope
one accounts for this: $P_{\Gm\times\Ga}(y)=\max(y-1,0)$ and
$P_{\Wh_2}(y)=y$. The first inequality is a form of upper semicontinuity of
fibre dimension. It is not a formal consequence of Chevalley's theorem,
because $G[p^m]$ is a formal scheme over $\Lambda$, not a scheme of finite
type.

The proof of Theorem~\ref{thm:A} (\S\ref{sec:geometric}) compares the finite
group schemes
\[
  H_{n,m}(G)=\ker\bigl(\Fr^n_{G/\Lambda}\bigr)\cap\ker\bigl([p^m]_G\bigr),
  \qquad m,n\ge1,
\]
on the two fibres. They are closed subschemes of the finite free Frobenius
kernels, cut out by finitely many equations, so their orders are upper
semicontinuous on $\Spec\Lambda$, even when $G[p^m]$ has positive dimension.
For fixed $m$ the order of $H_{n,m}$ grows like $p^{n\dim G[p^m]}$. Along
$m\approx(1-x)n$ its exponential growth rate is $\dim G-P_G^*(x)$, where
$P_G^*$ is the Legendre transform of $P_G$. Comparing Frobenius kernels with
kernels of multiplication by powers of $p$ is a familiar device in the theory
of $p$-divisible groups; compare \cite[Ch.~IV, \S7]{Dem72}.

\subsection{Cartier modules}
Let $A$ be a perfect $\Fp$-algebra and let $G$ be a smooth formal group over
$A$ with finite projective Lie algebra. Its covariant Cartier module
$M_A=\TC(G)$ of $p$-typical curves is a module over the Witt vectors $W(A)$,
with operators $F$ and~$V$. For a family over $k[[t]]$, with $k$ a perfect
field, one may try to prove Theorem~\ref{thm:A} by passing to the completed
perfection $R$ of $k[[t]]$ and comparing $M_R$ with the Cartier modules of the
two fibres by extension of Witt scalars. Part~(2) of the next theorem shows
that this fails for the full module, already for the constant family $\Ga$.
Part~(1) shows that it works for the quotients $M/V^nM$.

\begin{mainthm}\label{thm:B}
\begin{enumerate}[label=(\arabic*)]
\item Let $A\to B$ be a homomorphism of perfect $\Fp$-algebras, let $G$ be a
  smooth formal group over $A$ with finite projective Lie algebra, and put
  $M_A=\TC(G)$ and $M_B=\TC(G_B)$. For every $n\ge1$ the natural map
  \[
    (M_A/V^nM_A)\otimes_{W_n(A)}W_n(B)\longrightarrow M_B/V^nM_B
  \]
  is an isomorphism of finitely presented $W_n(B)$-modules. In particular,
  \[
    M_B\cong\varprojlim_n\,(M_A/V^nM_A)\otimes_{W_n(A)}W_n(B).
  \]
\item Let $k$ be a perfect field of characteristic $p$, let $R$ be the
  $t$-adic completion of the perfection of $k[[t]]$, and let $K$ be the
  completion of an algebraic closure of $k((t))$. For $G=\Ga$ over $R$, the
  natural map $M_R\otimes_{W(R)}W(K)\to M_K$ is not surjective, and the
  natural map $M_R\otimes_{W(R)}W(k)\to M_k$ is not injective.
\end{enumerate}
\end{mainthm}

Part~(1) is Theorem~\ref{thm:finite-level} and part~(2) is
Proposition~\ref{prop:Ga}. Appendix~\ref{app:lattice} shows that injectivity
on the special fibre also fails for the one-dimensional family $X+Y+tXY$,
whose generic fibre is $\Gm$. Part~(1) can be combined with a determinantal
argument parallel to \cite[Lemma~4.14]{BMS18}. For a smooth formal group over
$R$ this bounds the length of $M/(V^nM+p^mM)$ on the generic fibre by its
length on the special fibre, for every pair $(n,m)$ and with no error term
(Proposition~\ref{prop:Q-semicontinuity}). These lengths differ from
$\log_p|H_{n,m}|$ by a bounded amount (Corollary~\ref{cor:compare}). This
gives a second proof of Theorem~\ref{thm:A} for families over $k[[t]]$
(Corollary~\ref{cor:second-proof}).

\subsection{Explicit specializations}

\begin{mainthm}\label{thm:C}
\begin{enumerate}[label=(\arabic*)]
\item For every $r\ge1$ there is a smooth formal group over $\Fp[[t]]$ whose
  special fibre is isomorphic to $\Wh_{r+1}$ and whose geometric generic
  fibre is isomorphic to $\Gm\times\Wh_r$. For $r=1$ it may be given by a
  polynomial formal group law over $\Fp[t]$ whose generic fibre is
  isomorphic to $\Gm\times\Ga$ over $\Fp(t)$. For $p=2$ this law is
  \[
    (X,U)*(Y,Z)=\bigl(X+Y+tXY,\ U+Z+XY\bigr).
  \]
\item For all integers $a\ge b+2$ and $b\ge0$ there is a smooth formal group
  over $\Fp[[t]]$ whose special fibre is isomorphic to
  $\Wh_{a-1}\times\Wh_{b+1}$ and whose geometric generic fibre is isomorphic
  to $\Wh_a\times\Wh_b$, where $\Wh_0=0$.
\end{enumerate}
\end{mainthm}

These families are instances of the following converse to
Theorem~\ref{thm:A}.

\begin{question}\label{q:converse}
Let $G_1$ and $G_2$ be smooth commutative formal groups of the same dimension
over an algebraically closed field $\kappa$ of characteristic $p$, such that
$\dim G_1[p^m]\le\dim G_2[p^m]$ for all $m\ge1$ and $P_{G_1}\le P_{G_2}$. Is
there a smooth formal group over $\kappa[[t]]$ whose geometric generic fibre
is isogenous to $G_1$ and whose special fibre is isogenous to $G_2$?
\end{question}

For $p$-divisible groups, Grothendieck conjectured and Oort proved
\cite{Oor00} that every Newton polygon lying below that of a given
$p$-divisible group $X_0$, with the same endpoints, is the Newton polygon of
the generic fibre of some deformation of $X_0$. For unipotent groups, the
hypotheses of Question~\ref{q:converse} say that the partition of Witt
lengths of $G_2$ is dominated by that of $G_1$. The families of
Theorem~\ref{thm:C}(2), multiplied by constant groups, realize every covering
relation of the dominance order (Remark~\ref{rem:dominance}).

\subsection{Organization and notation}
Section~\ref{sec:prelim} collects conventions on Witt vectors and Cartier
modules, the standard groups $\Ga$, $\Gm$, $\Wh_r$ and $L_{d,h}$, and the
extended Newton polygon. Section~\ref{sec:geometric} proves
Theorem~\ref{thm:A}. Section~\ref{sec:basechange} proves Theorem~\ref{thm:B},
and Section~\ref{sec:lengths} proves the length inequalities and gives the
second proof of Theorem~\ref{thm:A}. Section~\ref{sec:families} proves
Theorem~\ref{thm:C}. Appendix~\ref{app:lattice} computes the Cartier module
of the family $X+Y+tXY$.

For $z\in\mathbf R$ we write $z_+=\max(z,0)$. Lengths of modules over a
discrete valuation ring $W(\kappa)$ are denoted $\length_{W(\kappa)}$. For a
finite group scheme $H$ over a field $\kappa$ we write
$|H|=\dim_\kappa\cO(H)$.

\section{Preliminaries}\label{sec:prelim}

\subsection{Witt vectors}
Let $A$ be a perfect $\Fp$-algebra. We write $W(A)$ for its ring of
$p$-typical Witt vectors, $W_n(A)$ for the Witt vectors of length $n$, $[a]$
for the Teichm\"uller representative of $a\in A$, and $\sigma$ for the Witt
vector Frobenius, an automorphism of $W(A)$. The ring $W(A)$ is
$p$-torsion-free and $p$-adically complete, and every element can be written
uniquely as $\sum_{i\ge0}p^i[a_i]$ with $a_i\in A$. The Verschiebung equals
$p\sigma^{-1}$, and $W_n(A)=W(A)/p^nW(A)$. If $A$ is a perfect field, then
$W(A)$ is a complete discrete valuation ring with uniformizer $p$ and residue
field $A$, and
\begin{equation}\label{eq:valuation}
  v_p\Bigl(\sum_{i\ge0}p^i[a_i]\Bigr)=\min\{i: a_i\neq0\}.
\end{equation}

\subsection{Formal groups and Cartier modules}
Let $A$ be a ring of characteristic $p$. By a \emph{smooth formal group} over
$A$ we mean a commutative formal Lie group over $A$ with finitely generated
projective Lie algebra, viewed as a functor on nilpotent $A$-algebras
\cite{Zin84}. If $\Lie(G)$ is free of rank $D$, for instance if $A$ is local,
then $G$ is given by a commutative formal group law
$\Phi(X,Y)\in A[[X,Y]]^D$ in $D$ variables $X=(X_1,\dots,X_D)$.

We use Cartier's theory of curves \cite{Car67}, in the form of
\cite[Ch.~IV]{Zin84}. A \emph{curve} in $G$ is an element of
$G(TA[[T]])=\varprojlim_NG(TA[T]/T^N)$. When $G$ is given by $\Phi$, a curve
is a $D$-tuple $\gamma(T)$ of power series without constant term, and curves
are added with $\Phi$. Cartier's operators on curves include
$[a]\gamma(T)=\gamma(aT)$ for $a\in A$, $V_\ell\gamma(T)=\gamma(T^\ell)$ and
the Frobenius operators $F_\ell$ for primes $\ell$. A curve is
\emph{$p$-typical} if $F_\ell\gamma=0$ for every prime $\ell\neq p$. We write
$\TC(G)$ for the group of $p$-typical curves. It is a left module over the
$p$-typical Cartier ring $\E_A$ of \cite[Ch.~IV]{Zin84}, which contains
$W(A)$ and elements $F=F_p$ and $V=V_p$ subject to
\begin{equation}\label{eq:cartier-relations}
  F\xi=\sigma(\xi)F,\qquad \xi V=V\sigma(\xi),\qquad
  V\xi F={}^{V}\xi,\qquad FV=p\qquad(\xi\in W(A)),
\end{equation}
where ${}^V\xi$ denotes the Witt vector Verschiebung. The Teichm\"uller
element $[a]\in W(A)$ acts on curves as above, and $V$ acts by
$\gamma(T)\mapsto\gamma(T^p)$. Since $A$ has characteristic $p$, we have
${}^V1=p$ in $W(A)$. Hence $VF=FV=p$ in $\E_A$, and therefore
\begin{equation}\label{eq:FV-commute}
  FV^j=V^jF,\qquad [a]V^j=V^j[a^{p^j}],\qquad {}^{V^j}\xi=V^j\xi F^j
  \qquad(j\ge0,\ a\in A,\ \xi\in W(A)).
\end{equation}

\begin{theorem}[Cartier; see {\cite[4.23, 4.39, 4.43, 4.45]{Zin84}}]
\label{thm:cartier}
Let $G$ be a smooth formal group over $A$ and let $M=\TC(G)$.
\begin{enumerate}[label=(\arabic*)]
\item The module $M$ is \emph{$V$-reduced}: $V$ is injective, the map
  $M\to\varprojlim_nM/V^nM$ is bijective, and the tangent map induces an
  isomorphism $M/VM\cong\Lie(G)$ of $A$-modules. The functor
  $G\mapsto\TC(G)$ is an equivalence between smooth formal groups over $A$
  and $V$-reduced $\E_A$-modules $M$ with $M/VM$ finitely generated
  projective over $A$.
\item Suppose that $\Lie(G)$ is free, and let $e_1,\dots,e_D\in M$ lift a
  basis of $M/VM$. Then every element of $M$ has a unique expansion
  $\sum_{i=1}^D\sum_{j\ge0}V^j[a_{ij}]e_i$ with $a_{ij}\in A$; we say that
  $(e_i)$ is a \emph{$V$-basis}. In particular, there are unique elements
  $c_{ijk}\in A$ with
  \begin{equation}\label{eq:structure}
    Fe_i=\sum_{j=1}^D\sum_{k\ge0}V^k[c_{ijk}]e_j\qquad(1\le i\le D).
  \end{equation}
  Conversely, for every family $(c_{ijk})$ of elements of $A$ there is a
  smooth formal group over $A$, unique up to isomorphism, whose Cartier
  module has a $V$-basis satisfying \eqref{eq:structure}.
\item For every homomorphism $A\to B$, base change of curves identifies
  $\TC(G_B)$ with $\E_B\otimes_{\E_A}M$. In the situation of (2), the images
  of $e_1,\dots,e_D$ form a $V$-basis of $\TC(G_B)$ that satisfies
  \eqref{eq:structure} with the images of the $c_{ijk}$ in $B$.
\end{enumerate}
\end{theorem}

The existence and uniqueness of the expansions in (2) follow from (1) by
successive approximation modulo $V^n$. The converse part of (2) is
\cite[4.39]{Zin84}: the module is the cokernel of the endomorphism of
$\bigoplus_i\E_Ae_i$ sending $e_i$ to $Fe_i-\sum_{j,k}V^k[c_{ijk}]e_j$.
Assertion~(3) follows from this presentation and \cite[4.43,
4.45]{Zin84}. No completion is needed in (3), because $M/VM$ is finitely
presented.

\begin{example}\label{ex:standard}
\begin{enumerate}[label=(\alph*), wide]
\item \emph{The additive group.} For the law $X+Y$, the operators $F_\ell$
  send $\sum_nc_nT^n$ to $\ell\sum_nc_{\ell n}T^n$. Hence the $p$-typical
  curves of $\Ga$ are the series $\sum_{i\ge0}a_iT^{p^i}$ with $a_i\in A$,
  and $F$ acts on them as zero. The curve $T$ is a $V$-basis with $FT=0$.
\item \emph{Witt groups.} For $r\ge1$ let $\Wh_r$ be the formal completion of
  the group scheme $W_r$ of Witt vectors of length $r$. It is a smooth formal
  group of dimension $r$ with coordinates $x_0,\dots,x_{r-1}$. Its Cartier
  module has a $V$-basis $e_1,\dots,e_r$ with $Fe_i=e_{i+1}$ for $i<r$ and
  $Fe_r=0$ \cite[4.40]{Zin84}. In characteristic $p$, multiplication by $p$
  on $W_r$ is $(x_0,\dots,x_{r-1})\mapsto(0,x_0^p,\dots,x_{r-2}^p)$. Hence
  \begin{equation}\label{eq:p-on-witt}
    [p^m](x_0,\dots,x_{r-1})
    =\bigl(0,\dots,0,x_0^{p^m},\dots,x_{r-m-1}^{p^m}\bigr),
  \end{equation}
  with $\min(m,r)$ zeros; in particular $[p^m]=0$ for $m\ge r$. We put
  $\Wh_0=0$.
\item \emph{The multiplicative group.} Let $\Gm$ have the law $X+Y+XY$, so
  that $1+X$ is the multiplicative coordinate, and let
  $E(T)=\exp\bigl(\sum_{i\ge0}T^{p^i}/p^i\bigr)\in1+T\Zpl[[T]]$ be the
  Artin--Hasse exponential. Over a $\mathbf Q$-algebra, the logarithm
  $\log(1+X)$ transforms $F_\ell$ into the map
  $\sum_nc_nT^n\mapsto\ell\sum_nc_{\ell n}T^n$. Applied to
  $\log E(T)=\sum_iT^{p^i}/p^i$, this shows that the curve
  $\varepsilon=E(T)-1$ is $p$-typical and satisfies $F\varepsilon=\varepsilon$.
  These identities have coefficients in $\Zpl$, so they hold over every $A$.
  The tangent vector of $\varepsilon$ is a basis of $\Lie(\Gm)$, so
  $\varepsilon$ is a $V$-basis. By \eqref{eq:FV-commute}, a Witt vector $w$
  with components $(w_0,w_1,\dots)$ acts by
  $w\varepsilon=\sum_iV^i[w_i]F^i\varepsilon=\sum_iV^i[w_i]\varepsilon$.
  Consequently $w\mapsto w\varepsilon$ is an isomorphism of $W(A)$-modules
  $W(A)\cong\TC(\Gm)$, given by
  \begin{equation}\label{eq:AH}
    w\longmapsto\prod_{i\ge0}E\bigl(w_iT^{p^i}\bigr)-1 .
  \end{equation}
  Under it, $F$ corresponds to $\sigma$ and $V$ to the Witt vector
  Verschiebung.
\end{enumerate}
\end{example}

\subsection{Isogenies and the classification}\label{subsec:isogeny}
Let $\kappa$ be a field of characteristic $p$. A homomorphism $u\colon G\to H$
of smooth formal groups of the same dimension over $\kappa$ is an
\emph{isogeny} if $\ker u$ is finite. In that case
$\cO(G)/\fm_H\cO(G)=\cO(\ker u)$ is finite dimensional. Hence $\cO(G)$ is a
finite $\cO(H)$-module, and it is free by the miracle flatness theorem,
because both rings are regular local rings of the same dimension. Its rank
is $|\ker u|=p^a$, and $p^a$ is the \emph{degree} of $u$. By
\cite[5.25]{Zin84} there are $j\ge0$ and $v\colon H\to G^{(p^j)}$ with
$v\circ u=\Fr^j_G$, where $\Fr_G\colon G\to G^{(p)}$ is the relative
Frobenius. We call two smooth formal groups \emph{isogenous} if they are
connected by a finite chain of isogenies.

\begin{construction}\label{con:Ldh}
Let $1\le d\le h$ be integers. Let $N_{d,h}$ be the free $\Zp$-module with
basis $e_0,\dots,e_{h-1}$, and define $e_i$ for all $i\ge0$ by
$e_{i+h}=pe_i$. Let $V$ and $F$ be the $\Zp$-linear endomorphisms with
$Ve_i=e_{i+d}$ and $Fe_i=e_{i+h-d}$. Then $FV=VF=p$ and $V^h=p^d$. In
particular, $V$ is injective and $N_{d,h}$ is $V$-adically complete, and
$N_{d,h}/VN_{d,h}$ has basis $e_0,\dots,e_{d-1}$. Since $W(\Fp)=\Zp$ and
$\sigma=\mathrm{id}$, these operators make $N_{d,h}$ a $V$-reduced
$\E_{\Fp}$-module. By Theorem~\ref{thm:cartier}(1) there is a smooth formal
group $L_{d,h}$ over $\Fp$ of dimension $d$ with $\TC(L_{d,h})=N_{d,h}$.
Comparing with Example~\ref{ex:standard}(c) shows that $L_{1,1}\cong\Gm$.
\end{construction}

\begin{lemma}\label{lem:Ldh}
Let $L=L_{d,h}$, and identify $L^{(p)}$ with $L$ using the $\Fp$-structure.
Then $\Fr_L^h=[p^d]_L$. Consequently, for every field $\kappa$ of
characteristic $p$, the group $L_\kappa$ is a $p$-divisible formal group of
height $h$ and dimension $d$, and $|L_\kappa[p^k]|=p^{hk}$ for all $k\ge0$.
\end{lemma}

\begin{proof}
For a curve $\gamma$ of $L$ over $\Fp$, the curve $\Fr_L\circ\gamma$ is
obtained by raising the coordinates of $\gamma(T)$ to the $p$-th power. Its
coefficients lie in $\Fp$, so $\Fr_L\circ\gamma=\gamma(T^p)=V\gamma$. Hence
$\TC(\Fr_L^h)=V^h=p^d=\TC([p^d])$, and $\Fr_L^h=[p^d]$ because $\TC$ is
faithful. After base change to $\kappa$ we get
$\ker[p^d]=\ker\Fr^h_{L_\kappa}$, which has order $p^{hd}$. Hence $[p]$ is an
isogeny, of degree $p^h$ because degrees are multiplicative, and
$|L_\kappa[p^k]|=p^{hk}$.
\end{proof}

\begin{theorem}[Classification up to isogeny]\label{thm:classification}
Let $\kappa$ be an algebraically closed field of characteristic $p$. Every
smooth formal group over $\kappa$ is isogenous to a product
\begin{equation}\label{eq:decomposition}
  \prod_{i=1}^u\Wh_{r_i}\times\prod_{j=1}^v(L_{d_j,h_j})_\kappa ,
  \qquad r_i\ge1,\quad 1\le d_j\le h_j,\quad \gcd(d_j,h_j)=1.
\end{equation}
\end{theorem}

\begin{proof}
By \cite[5.36]{Zin84}, every smooth formal group $G$ over $\kappa$ is an
extension of a $p$-divisible formal group by a unipotent one, and $G$ is
isogenous to their product. By \cite[5.38]{Zin84}, the unipotent factor is
isogenous to a product of groups $\Wh_r$. For the $p$-divisible factor $G'$,
the isocrystal $\TC(G')[1/p]$ is a direct sum of simple isocrystals by the
Dieudonn\'e--Manin theorem \cite{Man63}, \cite[Ch.~IV]{Dem72}. For coprime
$d$ and $h$, the isocrystal $N_{d,h}\otimes_{\Zp}W(\kappa)[1/p]$ is the simple
isocrystal on which $F$ has slope $(h-d)/h$. Isogeny classes of
$p$-divisible formal groups correspond to isomorphism classes of their
isocrystals, which gives the claim.
\end{proof}

\subsection{The invariants}\label{subsec:invariants}
Let $B$ be a smooth formal group of dimension $D$ over an algebraically
closed field $\kappa$, and choose an isogeny decomposition
\eqref{eq:decomposition}. We call $\rho(B)=(r_1,\dots,r_u)$, arranged in
decreasing order, the \emph{partition of Witt lengths} of $B$. We call
$\lambda_j=1-d_j/h_j\in[0,1)$ the \emph{slopes} of $B$, counted with
multiplicity $h_j$. On $N_{d,h}$ we have $F^h=p^{h-d}$, so $\lambda_j$ is the
slope of $F$ on the covariant Cartier module. It is one minus the slope in
the contravariant convention of \cite{Man63} and \cite{Dem72}. The group
$\Gm=L_{1,1}$ has slope~$0$, and the factors $\Wh_r$ contribute no slopes.

The \emph{extended Newton polygon} $P_B\colon[0,\infty)\to\mathbf R$ is the
continuous, convex, piecewise linear function with $P_B(0)=0$ whose slope is
$\lambda_j$ on an interval of length $h_j$ for each $j$, with the slopes in
increasing order, and equals $1$ on $[\sum_jh_j,\infty)$. A priori $\rho(B)$
and $P_B$ depend on the chosen decomposition. Corollary~\ref{cor:invariants}
shows that they depend only on the isogeny class of $B$.

\begin{lemma}\label{lem:legendre}
For $x\in[0,1]$ put $P_B^*(x)=\sup_{y\ge0}\bigl(xy-P_B(y)\bigr)$. Then
\[
  P_B^*(x)=\sum_jh_j(x-\lambda_j)_+\quad(0\le x\le1)
  \qquad\text{and}\qquad
  P_B(y)=\sup_{0\le x\le1}\bigl(xy-P_B^*(x)\bigr)\quad(y\ge0).
\]
Consequently, for two such groups $B$ and $B'$ we have $P_B\le P_{B'}$ on
$[0,\infty)$ if and only if $P_B^*\ge P_{B'}^*$ on $[0,1]$.
\end{lemma}

\begin{proof}
Fix $x\in[0,1]$. The function $y\mapsto xy-P_B(y)$ is concave and piecewise
linear. Its slope is $x-\lambda_j$ along the segments of slope $\lambda_j$,
and $x-1\le0$ along the final ray. Hence it attains its supremum at
$y_x=\sum_{\lambda_j<x}h_j$, where $P_B(y_x)=\sum_{\lambda_j<x}h_j\lambda_j$.
The supremum is therefore $\sum_jh_j(x-\lambda_j)_+$. For the second formula,
the inequality $P_B(y)\ge xy-P_B^*(x)$ holds by definition. Conversely, all
slopes of $P_B$ lie in $[0,1]$, so every $y\ge0$ has a supporting line of
some slope $x\in[0,1]$. For this $x$ we get $xz-P_B(z)\le xy-P_B(y)$ for all
$z\ge0$, i.e. $P_B^*(x)=xy-P_B(y)$. The last assertion follows from the two
formulas.
\end{proof}

\section{Specialization via finite Frobenius kernels}\label{sec:geometric}

\subsection{\texorpdfstring{The group schemes $H_{n,m}$}{The group schemes H(n,m)}}
Let $\Lambda$ be a ring of characteristic $p$, and let $G$ be a smooth formal
group over $\Lambda$ given by a formal group law $\Phi$ in $D$ variables
$X=(X_1,\dots,X_D)$. For $n\ge1$ let $G^{(p^n)}$ be given by the law obtained
from $\Phi$ by raising its coefficients to the $p^n$-th power, and let
$\Fr^n=\Fr^n_{G/\Lambda}\colon G\to G^{(p^n)}$, $X\mapsto X^{p^n}$, be the
relative Frobenius. Its kernel
\[
  \ker(\Fr^n)=\Spec\Lambda[X]/\bigl(X_1^{p^n},\dots,X_D^{p^n}\bigr)
\]
is a finite free group scheme of order $p^{nD}$. For $m\ge1$, let
$[p^m](X)=([p^m]_1(X),\dots,[p^m]_D(X))$ be the power series of
multiplication by $p^m$ on $G$. Since
$\Fr^n\circ[p^m]_G=[p^m]_{G^{(p^n)}}\circ\Fr^n$, the endomorphism $[p^m]_G$
preserves $\ker(\Fr^n)$. We put
\[
  H_{n,m}(G)=\ker\bigl([p^m]_G|_{\ker(\Fr^n)}\bigr)
           =\ker(\Fr^n)\cap\ker([p^m]_G).
\]
This is a finite group scheme over $\Lambda$ with coordinate ring
\begin{equation}\label{eq:Hnm}
  \cO(H_{n,m}(G))=\Lambda[X]/\bigl(X_1^{p^n},\dots,X_D^{p^n},\,
  [p^m]_1(X),\dots,[p^m]_D(X)\bigr),
\end{equation}
where $[p^m]_i(X)$ denotes the image of the power series in the truncated
polynomial ring $\Lambda[X]/(X_1^{p^n},\dots,X_D^{p^n})$. The formation of
$H_{n,m}(G)$
commutes with base change in $\Lambda$, and $H_{n,m}(G_1\times G_2)$ is the
product of $H_{n,m}(G_1)$ and $H_{n,m}(G_2)$. For a smooth formal group $B$
over a field $\kappa$ we put
\[
  q_{n,m}(B)=\log_p|H_{n,m}(B)|=\log_p\dim_\kappa\cO(H_{n,m}(B)).
\]
This number does not change under extension of the base field, and it is
additive on products.

\begin{lemma}\label{lem:semicontinuity}
Let $\Lambda$ be a local domain of characteristic $p$, let $G$ be as above,
and let $\bar\eta$ and $\bar s$ be geometric points over the generic point and
the closed point of $\Spec\Lambda$. Then
$q_{n,m}(G_{\bar\eta})\le q_{n,m}(G_{\bar s})$ for all $n,m\ge1$.
\end{lemma}

\begin{proof}
Let $P=\Lambda[X]/(X_1^{p^n},\dots,X_D^{p^n})$, a free $\Lambda$-module of
rank $p^{nD}$. By \eqref{eq:Hnm}, $\cO(H_{n,m}(G))$ is the cokernel of the
$\Lambda$-linear map $\theta\colon P^D\to P$,
$(g_1,\dots,g_D)\mapsto\sum_ig_i\,[p^m]_i(X)$. Let $\Theta$ be the matrix of
$\theta$ in monomial bases. For every homomorphism from $\Lambda$ to a field
$\kappa$, right exactness of the tensor product gives
$\dim_\kappa\cO(H_{n,m}(G_\kappa))=p^{nD}-\rk(\Theta\otimes\kappa)$. A minor
of $\Theta$ that is nonzero in the residue field of $\Lambda$ is nonzero in
$\Lambda$, hence in $\Frac(\Lambda)$, since $\Lambda$ is a domain. Thus the
rank of $\Theta$ at the closed point is at most its rank at the generic
point. Ranks do not change under field extensions.
\end{proof}

\begin{lemma}\label{lem:kernel-dimension}
Let $B$ be a smooth formal group over a field $\kappa$, given by a law in $D$
variables, and let $m\ge1$. Put
$C=\kappa[[X]]/([p^m]_1(X),\dots,[p^m]_D(X))$, so that $B[p^m]=\Spf C$, and
let $\dim B[p^m]$ be the Krull dimension of $C$. Then
\[
  \lim_{n\to\infty}\frac{q_{n,m}(B)}{n}=\dim B[p^m].
\]
\end{lemma}

\begin{proof}
We may assume $D\ge1$. Let $\fm$ be the maximal ideal of $C$; it is generated
by the images of $X_1,\dots,X_D$. Put $Q=p^n$. Since
$\kappa[[X]]/(X_1^Q,\dots,X_D^Q)=\kappa[X]/(X_1^Q,\dots,X_D^Q)$, formula
\eqref{eq:Hnm} gives $\cO(H_{n,m}(B))=C/(X_1^Q,\dots,X_D^Q)C$. Every monomial
of degree $D(Q-1)+1$ in $X_1,\dots,X_D$ is divisible by some $X_i^Q$, so
\[
  \fm^{DQ}\subseteq\fm^{D(Q-1)+1}\subseteq(X_1^Q,\dots,X_D^Q)C\subseteq\fm^Q,
\]
and therefore
$\length(C/\fm^Q)\le|H_{n,m}(B)|\le\length(C/\fm^{DQ})$. Let
$\delta=\dim C$. By the theory of Hilbert--Samuel polynomials there are
constants $0<c_1\le c_2$ with $c_1s^\delta\le\length(C/\fm^s)\le c_2s^\delta$
for all $s\ge1$. Hence $c_1p^{n\delta}\le|H_{n,m}(B)|\le c_2D^\delta
p^{n\delta}$. Taking $\log_p$ and dividing by $n$ proves the claim.
\end{proof}

\begin{lemma}\label{lem:isogeny-error}
Let $u\colon B\to B'$ be an isogeny of degree $p^a$ over a field $\kappa$.
Then $|q_{n,m}(B)-q_{n,m}(B')|\le a$ for all $n,m\ge1$.
\end{lemma}

\begin{proof}
Put $H=H_{n,m}(B)$ and $H'=H_{n,m}(B')$, and let $E=B\times_{B'}H'$ be the
scheme-theoretic preimage of $H'$. Since $\cO(B)$ is free of rank $p^a$ over
$\cO(B')$, the algebra $\cO(E)=\cO(B)\otimes_{\cO(B')}\cO(H')$ is free of rank
$p^a$ over $\cO(H')$, and $|E|=p^a|H'|$. For every point $x$ of $E$ we have
$u(p^mx)=p^mu(x)=0$ and $u^{(p^n)}(\Fr^nx)=\Fr^n(u(x))=0$. Hence
$f=([p^m],\Fr^n)$ is a homomorphism $E\to\ker u\times\ker u^{(p^n)}$. Its
kernel is $E\cap H=H$, because $u(H)\subseteq H'$. Orders of finite group
schemes over a field are multiplicative in exact sequences. Therefore
$|E|=|H|\cdot|f(E)|$ with $1\le|f(E)|\le|\ker u|\cdot|\ker u^{(p^n)}|=p^{2a}$.
Thus $q_{n,m}(B)-q_{n,m}(B')=a-\log_p|f(E)|$ lies in $[-a,a]$.
\end{proof}

\begin{proposition}\label{prop:blocks}
Let $\kappa$ be a field of characteristic $p$ and let $m,n\ge1$.
\begin{enumerate}[label=(\arabic*)]
\item For $r\ge1$ one has $q_{n,m}(\Wh_r)=rn-(n-m)_+(r-m)_+$.
\item For $1\le d\le h$ and $L=(L_{d,h})_\kappa$ one has
  $|q_{n,m}(L)-\min(nd,mh)|\le hd$.
\end{enumerate}
\end{proposition}

\begin{proof}
(1) By \eqref{eq:p-on-witt} and \eqref{eq:Hnm},
\[
  \cO(H_{n,m}(\Wh_r))=\kappa[x_0,\dots,x_{r-1}]
  \big/\bigl(x_i^{p^n}\ (0\le i<r),\ x_i^{p^m}\ (0\le i<r-m)\bigr),
\]
which has dimension $p^e$ with $e=n\min(m,r)+\min(n,m)(r-m)_+$. If $m\ge r$,
then $e=nr$. If $m<r$ and $n\ge m$, then $e=nm+m(r-m)=rn-(n-m)(r-m)$. If
$m<r$ and $n<m$, then $e=nm+n(r-m)=rn$. In all cases
$e=rn-(n-m)_+(r-m)_+$.

(2) By Lemma~\ref{lem:Ldh}, $\Fr^{bh}=[p^{bd}]$ on $L$ for $b\ge0$, and
$|L[p^k]|=p^{hk}$. Let $b=\lfloor n/h\rfloor$. Since
$\ker\Fr^{bh}\subseteq\ker\Fr^n\subseteq\ker\Fr^{(b+1)h}$, intersecting with
$L[p^m]$ gives
\[
  L\bigl[p^{\min(m,bd)}\bigr]\subseteq H_{n,m}(L)
  \subseteq L\bigl[p^{\min(m,(b+1)d)}\bigr],
\]
so that $h\min(m,bd)\le q_{n,m}(L)\le h\min(m,(b+1)d)$. Since
$nd-hd<hbd\le nd$, both bounds lie within $hd$ of $\min(mh,nd)$.
\end{proof}

\begin{corollary}\label{cor:invariants}
Let $\kappa$ be algebraically closed and let $B$ be a smooth formal group of
dimension $D$ over $\kappa$ with an isogeny decomposition
\eqref{eq:decomposition}. Put $\lambda_j=1-d_j/h_j$.
\begin{enumerate}[label=(\arabic*)]
\item There is a constant $c$, independent of $n$ and $m$, such that
  \[
    \Bigl|\,q_{n,m}(B)-\sum_{i}\bigl(r_in-(n-m)_+(r_i-m)_+\bigr)
    -\sum_j\min(nd_j,mh_j)\Bigr|\le c\qquad(n,m\ge1).
  \]
\item $\dim B[p^m]=\sum_i\min(m,r_i)$ for every $m\ge1$.
\item For $0\le x<1$ and $m_n=\lfloor(1-x)n\rfloor$,
  \[
    \lim_{n\to\infty}\frac{q_{n,m_n}(B)}{n}=D-\sum_jh_j(x-\lambda_j)_+
    =D-P_B^*(x).
  \]
\item The partition $\rho(B)$ and the polygon $P_B$ do not depend on the
  chosen decomposition. They, and the dimensions $\dim B[p^m]$, depend only
  on the isogeny class of $B$.
\end{enumerate}
\end{corollary}

\begin{proof}
(1) follows from the additivity of $q_{n,m}$, from
Lemma~\ref{lem:isogeny-error} applied along a chain of isogenies, and from
Proposition~\ref{prop:blocks}. For fixed $m$, (1) gives
$q_{n,m}(B)/n\to\sum_i\bigl(r_i-(r_i-m)_+\bigr)=\sum_i\min(m,r_i)$, so (2)
follows from Lemma~\ref{lem:kernel-dimension}. For (3), take $n$ so large
that $m_n\ge\max_ir_i$. Then the $i$-th summand in (1) equals $r_in$, and
$\min(nd_j,m_nh_j)/n\to\min\bigl(d_j,(1-x)h_j\bigr)=d_j-h_j(x-\lambda_j)_+$.
Since $D=\sum_ir_i+\sum_jd_j$, the first equality follows; the second is
Lemma~\ref{lem:legendre}.

For (4), note that $q_{n,m}(B)$ and $\dim B[p^m]$ are defined without a
decomposition, and that by Lemma~\ref{lem:isogeny-error} the function
$q_{n,m}$ changes by a bounded amount under isogeny. The partition $\rho(B)$
is determined by the function $m\mapsto\sum_i\min(m,r_i)=\dim B[p^m]$,
because $\#\{i:r_i\ge m\}=\dim B[p^m]-\dim B[p^{m-1}]$, with
$\dim B[p^0]=0$. By (3), $q_{n,m}(B)$ determines $P_B^*$ on $[0,1)$. By
continuity it determines $P_B^*$ on $[0,1]$, and hence $P_B$ by
Lemma~\ref{lem:legendre}.
\end{proof}

\begin{proof}[Proof of Theorem~\ref{thm:A}]
Both geometric fibres have dimension $D$. By Lemma~\ref{lem:semicontinuity},
$q_{n,m}(G_{\bar\eta})\le q_{n,m}(G_{\bar s})$ for all $n,m\ge1$. Dividing by
$n$ and letting $n\to\infty$, Lemma~\ref{lem:kernel-dimension} gives
$\dim G_{\bar\eta}[p^m]\le\dim G_{\bar s}[p^m]$. Applying the same
inequality with $m=m_n$, Corollary~\ref{cor:invariants}(3) gives
$P^*_{G_{\bar s}}(x)\le P^*_{G_{\bar\eta}}(x)$ for $0\le x<1$. By continuity
this also holds for $x=1$, and Lemma~\ref{lem:legendre} gives
$P_{G_{\bar\eta}}\le P_{G_{\bar s}}$.
\end{proof}

\section{Base change of Cartier modules}\label{sec:basechange}

Let $G$ be a smooth formal group over a perfect $\Fp$-algebra $A$, and let
$A\to B$ be a homomorphism of perfect $\Fp$-algebras. We write $M_A=\TC(G)$
and $M_B=\TC(G_B)$. Base change of curves is a $W(A)$-linear map $M_A\to M_B$
commuting with $F$ and $V$. It induces $W(B)$-linear maps
$M_A\otimes_{W(A)}W(B)\to M_B$ and
$(M_A/V^nM_A)\otimes_{W(A)}W(B)\to M_B/V^nM_B$, which we call the natural
maps.

Throughout this section, Section~\ref{sec:lengths} and
Appendix~\ref{app:lattice}, $k$ is a perfect field of characteristic $p$ and
$R$ is the $t$-adic completion of the perfection
$\bigcup_{n}k[[t^{1/p^n}]]$ of $k[[t]]$. It is a perfect valuation ring of
rank one with value group $\mathbf Z[1/p]$, maximal ideal $\fm$ and residue
field $k$. We normalize its valuation $v$ by $v(t)=1$. Let $K$ be the
completion of an algebraic closure of $k((t))$. It is algebraically closed,
and the perfection of $k[[t]]$, together with its $t$-adic completion $R$,
embeds into $K$. We fix such an embedding $R\subseteq K$.

\subsection{The full Cartier module}

\begin{proposition}\label{prop:Ga}
For $G=\Ga$ over $R$, the natural map $M_R\otimes_{W(R)}W(K)\to M_K$ is
injective but not surjective, and the natural map
$M_R\otimes_{W(R)}W(k)\to M_k$ is surjective but not injective.
\end{proposition}

\begin{proof}
Let $A$ be a perfect $\Fp$-algebra. By Example~\ref{ex:standard}(a), the
elements of $M_A=\TC(\Ga)$ are the series $f(T)=\sum_{i\ge0}a_iT^{p^i}$. As
$\Ga$ is killed by $p$, so is $M_A$, and $W(A)$ acts through $W(A)/p=A$. The
element $c\in A$ acts as $[c]$, that is, by $f(T)\mapsto f(cT)$. Hence the
map
\[
  \beta_A\colon M_A\longrightarrow\prod_{i\ge0}A,\qquad
  \sum_{i\ge0}a_iT^{p^i}\longmapsto\bigl(a_i^{1/p^i}\bigr)_{i\ge0},
\]
is an isomorphism of $A$-modules, natural in the perfect ring $A$. Since $M_R$
is killed by $p$, we have $M_R\otimes_{W(R)}W(B)=M_R\otimes_RB$ for
$B\in\{K,k\}$. Under $\beta$, the natural maps become the maps
$\alpha_B\colon(\prod_iR)\otimes_RB\to\prod_iB$.

For $B=K$, the source of $\alpha_K$ is
$(\prod_iR)[1/t]\otimes_{\Frac(R)}K$. Choosing a $\Frac(R)$-basis of $K$
shows that $\alpha_K$ is injective. Every element of the source is a finite
sum $\sum_\nu x_\nu\otimes c_\nu$ with $x_\nu\in\prod_iR$ and $c_\nu\in K$.
Its image $\sum_\nu c_\nu x_\nu$ has all coordinates of valuation at least
$\min_\nu v(c_\nu)$. The vector $(t^{-i})_{i\ge0}$ has no such bound. It
corresponds to the curve $\sum_{i\ge0}t^{-ip^i}T^{p^i}\in M_K$, which is
therefore not in the image.

For $B=k$, the source of $\alpha_k$ is $(\prod_iR)/\fm(\prod_iR)$, and
$\alpha_k$ reduces each coordinate modulo $\fm$. So $\alpha_k$ is surjective.
Let $b=(t^{1/p^i})_{i\ge0}$. Then $\alpha_k(\bar b)=0$. On the other hand, an
element of $\fm(\prod_iR)$ is a finite sum $\sum_\nu c_\nu x_\nu$ with
$c_\nu\in\fm$ and $x_\nu\in\prod_iR$. All of its coordinates have valuation
at least $\min_\nu v(c_\nu)>0$. Since $v(t^{1/p^i})=p^{-i}\to0$, we get
$b\notin\fm(\prod_iR)$. In terms of curves, the class of
$t\sum_{i\ge0}T^{p^i}$ in $M_R\otimes_{W(R)}W(k)$ is nonzero, but its image
in $M_k$ is zero.
\end{proof}

\subsection{Truncations}

\begin{theorem}\label{thm:finite-level}
Let $A\to B$ be a homomorphism of perfect $\Fp$-algebras, and let $G$ be a
smooth formal group over $A$. For every $n\ge1$ the natural map
\[
  \phi_n\colon(M_A/V^nM_A)\otimes_{W_n(A)}W_n(B)\longrightarrow M_B/V^nM_B
\]
is an isomorphism, and $M_A/V^nM_A$ is a finitely presented
$W_n(A)$-module. Consequently
$M_B\cong\varprojlim_n(M_A/V^nM_A)\otimes_{W_n(A)}W_n(B)$.
\end{theorem}

\begin{proof}
Put $P=\Lie(G)$. Then $\Lie(G_B)=P\otimes_AB$, and $M_A/VM_A=P$ and
$M_B/VM_B=P\otimes_AB$ by Theorem~\ref{thm:cartier}(1).

\emph{Step 1: graded pieces.} Let $i\ge0$. Since $V$ is injective, $V^i$
induces a bijection $M_A/VM_A\to V^iM_A/V^{i+1}M_A$. The quotient
$V^iM_A/V^{i+1}M_A$ is killed by $p$, because $pV^i=V^{i+1}F$. Hence $W(A)$
acts on it through $W(A)/p=A$, and by \eqref{eq:FV-commute} the element
$c\in A$ acts through $[c]V^i=V^i[c^{p^i}]$. Thus $V^i$ is an isomorphism of
$A$-modules
\[
  (\sigma^i)_*P\xrightarrow{\ \sim\ }V^iM_A/V^{i+1}M_A ,
\]
where $(\sigma^i)_*P$ denotes $P$ with $c\in A$ acting by multiplication with
$c^{p^i}$. Since $\sigma$ is an automorphism of $A$, the module
$(\sigma^i)_*P$ is finitely generated and projective. The same holds over
$B$. On these graded pieces, the natural map is
\[
  (\sigma^i)_*P\otimes_AB\longrightarrow(\sigma^i)_*(P\otimes_AB),\qquad
  u\otimes b\longmapsto u\otimes b^{p^i}.
\]
It is an isomorphism. For $P=A$ it becomes the identity of $B$ under the
isomorphisms $u\otimes b\mapsto u^{1/p^i}b$ and $w\mapsto w^{1/p^i}$, and
the general case follows by additivity.

\emph{Step 2: Tor vanishing.} Since $W(A)$ is $p$-torsion-free, the sequence
$0\to W(A)\xrightarrow{p}W(A)\to A\to0$ is a free resolution of $A$. As
$W(B)$ is also $p$-torsion-free, $\Tor_1^{W(A)}(A,W(B))=W(B)[p]=0$, and the
higher Tor groups vanish. The same holds for every direct summand of a finite
free $A$-module, in particular for the graded pieces of Step~1.

\emph{Step 3: induction.} Consider the exact sequence
\[
  0\to V^nM_A/V^{n+1}M_A\to M_A/V^{n+1}M_A\to M_A/V^nM_A\to0 .
\]
By Step~2 and induction on $n$, $\Tor_1^{W(A)}(M_A/V^nM_A,W(B))=0$ for all
$n$. Hence the sequence stays exact after tensoring with $W(B)$. The natural
maps form a morphism to the analogous exact sequence over $B$. By Step~1 and
the five lemma, induction on $n$ shows that the natural map
$(M_A/V^nM_A)\otimes_{W(A)}W(B)\to M_B/V^nM_B$ is an isomorphism for every
$n$. Since $p^n=V^nF^n$ kills $M_A/V^nM_A$ and $W(B)/p^n=W_n(B)$, this map
is $\phi_n$.

\emph{Step 4: finite presentation.} The $W_n(A)$-module $A=W_n(A)/pW_n(A)$
is finitely presented. Hence so is every direct summand of a finite free
$A$-module. Extensions of finitely presented modules are finitely presented,
so Step~1 shows that $M_A/V^nM_A$ is finitely presented over $W_n(A)$.

Finally, $M_B$ is $V$-adically separated and complete by
Theorem~\ref{thm:cartier}(1), which gives the last assertion.
\end{proof}

\begin{corollary}\label{cor:naive}
In the situation of Theorem~\ref{thm:finite-level}, let
$T=W(B)\otimes_{W(A)}M_A$, and put $V(b\otimes u)=\sigma^{-1}(b)\otimes Vu$.
This $V$ is well defined, and
$T/V^nT\cong(M_A/V^nM_A)\otimes_{W_n(A)}W_n(B)$. The natural map $T\to M_B$
identifies $M_B$ with the $V$-adic completion $\varprojlim_nT/V^nT$. In
particular, $T\to M_B$ is injective if and only if $\bigcap_nV^nT=0$, and it
is bijective if and only if $T$ is $V$-adically separated and complete.
\end{corollary}

\begin{proof}
The map $V$ is well defined because $\xi V=V\sigma(\xi)$ for
$\xi\in W(A)$. Since $\sigma$ is bijective on $W(B)$, $V^nT$ is the image of
$W(B)\otimes_{W(A)}V^nM_A$. So $T/V^nT=W(B)\otimes_{W(A)}(M_A/V^nM_A)$ by
right exactness, and Theorem~\ref{thm:finite-level} applies.
\end{proof}

\begin{remark}\label{rem:zink-base-change}
By Theorem~\ref{thm:cartier}(3), $M_B$ is the base change
$\E_B\otimes_{\E_A}M_A$ along the Cartier rings, and no completion is needed
there. The ring $\E_B$ contains the $V$-adically convergent series
$\sum_rV^r[b_r]$. Corollary~\ref{cor:naive} says that for perfect rings, the
Witt vector base change $T$ differs from $M_B$ only by $V$-adic completion.
Proposition~\ref{prop:Ga} shows that $T$ can fail to be complete (for $B=K$)
and fail to be separated (for $B=k$).
\end{remark}

\section{Lengths of truncated Cartier modules}\label{sec:lengths}

For a smooth formal group $B$ of dimension $D$ over a perfect field $\kappa$,
with $M=\TC(B)$, we put
\[
  \ell_{n,m}(B)=\length_{W(\kappa)}\bigl(M/(V^nM+p^mM)\bigr)\qquad(n,m\ge1).
\]
This is finite: by Step~1 of the proof of Theorem~\ref{thm:finite-level},
$M/V^nM$ has a filtration with $n$ graded pieces of dimension $D$ over
$\kappa$, so $\length_{W(\kappa)}(M/V^nM)=nD$.

\subsection{Semicontinuity}

\begin{lemma}\label{lem:determinantal}
Let $R_0$ be a perfect $\Fp$-algebra, let $\iota\colon R_0\to K_0$ be an
injective homomorphism to a perfect field, and let $\pi\colon R_0\to k_0$ be
a homomorphism to a perfect field. Let $Q$ be a finitely presented
$W_n(R_0)$-module, and put $Q_{K_0}=Q\otimes_{W_n(R_0)}W_n(K_0)$ and
$Q_{k_0}=Q\otimes_{W_n(R_0)}W_n(k_0)$. Then
\[
  \length_{W(K_0)}(Q_{K_0})\le\length_{W(k_0)}(Q_{k_0}).
\]
\end{lemma}

\begin{proof}
Choose a presentation $W_n(R_0)^a\xrightarrow{\Phi}W_n(R_0)^b\to Q\to0$,
lift $\Phi$ to a matrix $\tilde\Phi$ over $W(R_0)$, and let
$\Psi=(\tilde\Phi\mid p^nI_b)$, a $b\times(a+b)$ matrix. Then $Q$ is the
cokernel of $\Psi\colon W(R_0)^{a+b}\to W(R_0)^b$. Since $p^nQ=0$ and
$W(F)/p^n=W_n(F)$ for $F\in\{K_0,k_0\}$, the module $Q_F$ is the cokernel of
the image $\Psi_F$ of $\Psi$ over the discrete valuation ring $W(F)$. The
matrix $\Psi_F$ has rank $b$, because it contains $p^nI_b$. By the theory of
elementary divisors, $\length_{W(F)}(Q_F)$ is the minimum of
$v_p(\delta_F)$ over all $b\times b$ minors $\delta$ of $\Psi$. Write a minor
as $\delta=\sum_ip^i[\delta_i]$ with $\delta_i\in R_0$. By
\eqref{eq:valuation} and the injectivity of $\iota$,
\[
  v_p(\delta_{K_0})=\min\{i:\iota(\delta_i)\neq0\}=\min\{i:\delta_i\neq0\}
  \le\min\{i:\pi(\delta_i)\neq0\}=v_p(\delta_{k_0}).
\]
Taking the minimum over all minors proves the lemma.
\end{proof}

For $R_0=\mathcal O_{C^\flat}$, the tilt of the ring of integers of a
complete algebraically closed field $C$, this inequality is
\cite[Lemma~4.14]{BMS18}. There it is proved by d\'evissage, together with
an identity involving $\Tor_1$.

\begin{proposition}\label{prop:Q-semicontinuity}
Let $G$ be a smooth formal group over $R$, and let $n,m\ge1$. The module
$Q_{n,m}(G)=M_R/(V^nM_R+p^mM_R)$ is a finitely presented $W_n(R)$-module,
and
\begin{align*}
  Q_{n,m}(G)\otimes_{W_n(R)}W_n(K)&\cong M_K/(V^nM_K+p^mM_K),\\
  Q_{n,m}(G)\otimes_{W_n(R)}W_n(k)&\cong M_k/(V^nM_k+p^mM_k).
\end{align*}
In particular $\ell_{n,m}(G_K)\le\ell_{n,m}(G_k)$.
\end{proposition}

\begin{proof}
We have $Q_{n,m}(G)=(M_R/V^nM_R)/p^m(M_R/V^nM_R)$, which is finitely
presented by Theorem~\ref{thm:finite-level}. The base change formulas follow
from Theorem~\ref{thm:finite-level} and the right exactness of the tensor
product. The inequality follows from Lemma~\ref{lem:determinantal} applied
to the inclusion $R\subseteq K$ and the residue map $R\to k$.
\end{proof}

\subsection{Asymptotics}

\begin{proposition}\label{prop:length-blocks}
Let $\kappa$ be a perfect field of characteristic $p$, and let $n,m\ge1$.
\begin{enumerate}[label=(\arabic*)]
\item If $\kappa'\supseteq\kappa$ is a perfect field, then
  $\ell_{n,m}(B_{\kappa'})=\ell_{n,m}(B)$.
\item Let $u\colon B\to B'$ be an isogeny. Then $\TC(u)$ is injective, and
  its cokernel has finite length $c$. Moreover,
  $|\ell_{n,m}(B)-\ell_{n,m}(B')|\le c$.
\item $\ell_{n,m}(\Wh_r)=rn-(n-m)_+(r-m)_+$ for $r\ge1$.
\item $\ell_{n,m}((L_{d,h})_\kappa)=\min(nd,mh)$ for $1\le d\le h$.
\item $\ell_{n,m}(B_1\times B_2)=\ell_{n,m}(B_1)+\ell_{n,m}(B_2)$.
\end{enumerate}
\end{proposition}

\begin{proof}
(1) Let $Q=\TC(B)/(V^n+p^m)$. By Theorem~\ref{thm:finite-level} and right
exactness, $\TC(B_{\kappa'})/(V^n+p^m)\cong Q\otimes_{W(\kappa)}W(\kappa')$.
The ring $W(\kappa')$ is flat over the discrete valuation ring $W(\kappa)$,
being torsion-free, and $W(\kappa')/p=\kappa'$. By d\'evissage,
$\length_{W(\kappa')}(Q\otimes_{W(\kappa)}W(\kappa'))=\length_{W(\kappa)}(Q)$.

(2) Write $M=\TC(B)$ and $N=\TC(B')$. A curve $\gamma$ in $B$ corresponds to
a continuous homomorphism $\gamma^*\colon\cO(B)\to\kappa[[T]]$ mapping the
maximal ideal $\fm_B$ into $T\kappa[[T]]$. If $u\circ\gamma=0$, then
$\gamma^*$ factors through the artinian local ring $\cO(\ker u)$. Then
$\gamma^*(\fm_B)$ consists of nilpotent elements of $\kappa[[T]]$, so
$\gamma=0$. Thus $\TC(u)$ is injective, and we regard $M$ as a submodule of
$N$.

Next choose $v\colon B'\to B^{(p^j)}$ with $v\circ u=\Fr^j_B$
(\S\ref{subsec:isogeny}). Then
$u^{(p^j)}\circ v\circ u=u^{(p^j)}\circ\Fr^j_B=\Fr^j_{B'}\circ u$. Since $u$
is faithfully flat, it is an epimorphism, and so
$u^{(p^j)}\circ v=\Fr^j_{B'}$. For every $B''$ over $\kappa$ we have
$\TC(\Fr^j_{B''})(\gamma)=V^j\gamma^{(p^j)}$, where
$\gamma\mapsto\gamma^{(p^j)}$ raises the coefficients of $\gamma$ to the
$p^j$-th power. Since $\kappa$ is perfect, this map is a semilinear
bijection $\TC(B'')\to\TC(B''^{(p^j)})$, and it commutes with $V$. Hence the
image of $\TC(u^{(p^j)})$ contains $V^j\TC(B'^{(p^j)})$. Since
$(u\circ\gamma)^{(p^j)}=u^{(p^j)}\circ\gamma^{(p^j)}$, it follows that
$M\supseteq V^jN$. Therefore $c=\length(N/M)\le\length(N/V^jN)$ is finite.

Put $Q(X)=V^nX+p^mX$ for $X\in\{M,N\}$. Then $Q(M)\subseteq Q(N)$, and
computing $\length(N/Q(M))$ in two ways gives
\[
  \ell_{n,m}(B)-\ell_{n,m}(B')=\length\bigl(Q(N)/Q(M)\bigr)-c .
\]
The maps $x\mapsto V^nx$ and $x\mapsto p^mx$ from $N$ to $Q(N)/Q(M)$ vanish
on $M$, and their images generate $Q(N)/Q(M)$. The map $x\mapsto p^mx$ is
$W(\kappa)$-linear, so its image is a quotient of $N/M$. The map
$x\mapsto V^nx$ is $\sigma^{-n}$-semilinear, so its image is a quotient of
the twist of $N/M$ by the automorphism $\sigma^{-n}$ of $W(\kappa)$, which
also has length $c$. Hence $0\le\length(Q(N)/Q(M))\le2c$, and (2) follows.

(3) Let $e_1,\dots,e_r$ be the $V$-basis of Example~\ref{ex:standard}(b),
and let $M=\TC(\Wh_r)$. For $S\subseteq\{1,\dots,r\}\times\mathbf Z_{\ge0}$,
let $M_S$ be the set of elements whose expansion
$\sum_{i,j}V^j[a_{ij}]e_i$ has $a_{ij}=0$ for $(i,j)\notin S$. By
\eqref{eq:FV-commute}, $F(V^j[a]e_i)=V^j[a^p]Fe_i$. Since $F$ is additive and
continuous, and $\kappa$ is perfect, we get $F^mM=M_{\{i>m\}}=:M'$. The
module $M'$ is $V$-reduced with $V$-basis $e_{m+1},\dots,e_r$, and the
uniqueness of expansions gives $M'\cap V^sM=V^sM'$ for $s\ge0$. Now
$p^mM=V^mF^mM=V^mM'$ and $\length(M/V^nM)=nr$. If $n\le m$, then
$p^mM\subseteq V^nM$, and $\ell_{n,m}(\Wh_r)=nr$. If $n>m$, then
$p^mM\cap V^nM=V^m(M'\cap V^{n-m}M)=V^{n}M'$, because $V$ is injective.
Since $V^m$ is a semilinear bijection, it follows that
\[
  \ell_{n,m}(\Wh_r)=nr-\length\bigl(V^mM'/V^nM'\bigr)
  =nr-\length\bigl(M'/V^{n-m}M'\bigr)=nr-(n-m)(r-m)_+ .
\]

(4) Let $L=(L_{d,h})_\kappa$. By Corollary~\ref{cor:naive}, applied to
$\Fp\to\kappa$, the module $\TC(L)$ is the $V$-adic completion of
$W(\kappa)\otimes_{\Zp}N_{d,h}$, with $V(a\otimes e)=\sigma^{-1}(a)\otimes
Ve$. This module is free of finite rank over $W(\kappa)$ and satisfies
$V^h=p^d\sigma^{-h}$, so it is already $V$-adically separated and complete,
and it equals $\TC(L)$. For $s\ge0$ let
$T_{\ge s}=\sum_{j\ge s}W(\kappa)e_j$, which is free with basis
$e_s,\dots,e_{s+h-1}$. These submodules are nested, and
$T_{\ge s}/T_{\ge s+1}\cong W(\kappa)/p$, because $e_{s+h}=pe_s$. Moreover
$V^n\TC(L)=T_{\ge nd}$ and $p^m\TC(L)=T_{\ge mh}$. Hence
$V^n\TC(L)+p^m\TC(L)=T_{\ge\min(nd,mh)}$, which has colength
$\min(nd,mh)$.

(5) holds because $\TC(B_1\times B_2)=\TC(B_1)\oplus\TC(B_2)$.
\end{proof}

\begin{corollary}\label{cor:compare}
Let $B$ be a smooth formal group of dimension $D$ over an algebraically
closed field of characteristic $p$. Then $|\ell_{n,m}(B)-q_{n,m}(B)|$ is
bounded independently of $n,m\ge1$. Consequently
$\lim_n\ell_{n,m}(B)/n=\dim B[p^m]$ for $m\ge1$, and
$\lim_n\ell_{n,\lfloor(1-x)n\rfloor}(B)/n=D-P_B^*(x)$ for $0\le x<1$.
\end{corollary}

\begin{proof}
By Theorem~\ref{thm:classification} and
Proposition~\ref{prop:length-blocks}(2)--(5), $\ell_{n,m}(B)$ lies within
bounded distance of the explicit function in
Corollary~\ref{cor:invariants}(1). So does $q_{n,m}(B)$, by that corollary.
The limits follow from Corollary~\ref{cor:invariants}(2), (3) and
Lemma~\ref{lem:kernel-dimension}.
\end{proof}

\begin{corollary}\label{cor:second-proof}
Let $G_0$ be a smooth formal group over $k[[t]]$. For every pair
$n,m\ge1$,
\[
  \ell_{n,m}(G_{0,\bar\eta})\le\ell_{n,m}(G_{0,\bar s}).
\]
Together with Corollary~\ref{cor:compare}, this gives a second proof of
Theorem~\ref{thm:A} for $\Lambda=k[[t]]$, in which
Proposition~\ref{prop:Q-semicontinuity} replaces
Lemma~\ref{lem:semicontinuity}.
\end{corollary}

\begin{proof}
Let $G=G_0\otimes_{k[[t]]}R$. Its fibre over $k$ is the special fibre of
$G_0$. Its fibre over $K$ is $G_{0,\bar\eta}\otimes K$, where $\bar\eta$ is
given by the algebraic closure of $k((t))$ in $K$. By
Proposition~\ref{prop:Q-semicontinuity} and
Proposition~\ref{prop:length-blocks}(1),
\[
  \ell_{n,m}(G_{0,\bar\eta})=\ell_{n,m}(G_K)\le\ell_{n,m}(G_k)
  =\ell_{n,m}(G_{0,\bar s}).
\]
Dividing by $n$ and passing to the limits of Corollary~\ref{cor:compare}, the
argument of the proof of Theorem~\ref{thm:A} applies verbatim.
\end{proof}

\section{Explicit specializations}\label{sec:families}

\subsection{A polynomial formal group law}\label{subsec:polynomial-law}
Let
\[
  L(U)=\sum_{j=1}^{p-1}\frac{(-1)^{j-1}U^j}{j}\in\Zpl[U],\qquad
  S(X,Y)=X+Y+tXY,
\]
and put
\[
  C_t(X,Y)=\frac{L(tS)-L(tX)-L(tY)}{t^p},\qquad H_t(X)=\frac{L(tX)}{t^p}.
\]

\begin{lemma}\label{lem:Ct}
\begin{enumerate}[label=(\arabic*)]
\item $C_t\in\Zpl[t,X,Y]$. It is symmetric in $X$ and $Y$, satisfies
  $C_t(X,0)=0$, and all its monomials have degree at least $2$ in $(X,Y)$.
\item $C_t(X,Y)\equiv(-1)^p\,\dfrac{(X+Y)^p-X^p-Y^p}{p}\pmod{t}$, where the
  fraction is computed in $\mathbf Z[X,Y]$.
\item $C_t(X,Y)+C_t\bigl(S(X,Y),W\bigr)=C_t(Y,W)+C_t\bigl(X,S(Y,W)\bigr)$.
\end{enumerate}
\end{lemma}

\begin{proof}
(1) Since $1+tS=(1+tX)(1+tY)$, we have
$\log(1+tS)=\log(1+tX)+\log(1+tY)$ in $\mathbf Q[X,Y][[t]]$. As
$\log(1+U)-L(U)=\sum_{j\ge p}(-1)^{j-1}U^j/j$, we get
\begin{equation}\label{eq:numerator}
  L(tS)-L(tX)-L(tY)
  =-\sum_{j\ge p}\frac{(-1)^{j-1}}{j}\bigl((tS)^j-(tX)^j-(tY)^j\bigr).
\end{equation}
The right side is divisible by $t^p$, and the left side lies in
$\Zpl[t,X,Y]$. Hence $C_t\in\Zpl[t,X,Y]$. Symmetry and $C_t(X,0)=0$ are
clear. For $j\ge1$, the polynomial $S^j-X^j-Y^j$ has no monomials of degree
$<2$ in $(X,Y)$, which gives the last claim.

(2) In \eqref{eq:numerator}, the coefficient of $t^p$ comes from $j=p$ only,
and equals $-\frac{(-1)^{p-1}}{p}\bigl((X+Y)^p-X^p-Y^p\bigr)$.

(3) In $\Zpl[t^{\pm1},X,Y]$ we have
$C_t=H_t(S(X,Y))-H_t(X)-H_t(Y)$. For such a coboundary, the asserted
identity follows formally from the associativity of $S$. It holds in
$\Zpl[t,X,Y,W]$ because this ring injects into its localization at $t$.
\end{proof}

The identities of Lemma~\ref{lem:Ct} were also verified by computer for
$p\le7$. For instance, $C_t=XY$ for $p=2$, and
$C_t=-X^2Y-XY^2-\tfrac12tX^2Y^2$ for $p=3$.

\begin{proposition}\label{prop:polynomial-law}
Let $\bar C_t\in\Fp[t,X,Y]$ and $\bar H_t\in\Fp[t^{\pm1},X]$ be the reductions
of $C_t$ and $H_t$. Then
\[
  (X,U)*(Y,Z)=\bigl(X+Y+tXY,\ U+Z+\bar C_t(X,Y)\bigr)
\]
is a two-dimensional commutative formal group law over $\Fp[t]$, with the
following properties.
\begin{enumerate}[label=(\arabic*)]
\item Over $\Fp[t^{\pm1}]$, the change of coordinates
  $\psi(X,U)=\bigl(tX,\ U-\bar H_t(X)\bigr)$ is an isomorphism onto
  $\Gm\times\Ga$, with laws $x+y+xy$ and $u+z$.
\item At $t=0$ the law is the addition law of $W_2$ in Witt coordinates, so
  the special fibre is $\Wh_2$.
\item For $p=2$, $\bar C_t=XY$.
\end{enumerate}
\end{proposition}

\begin{proof}
The first coordinate is the commutative formal group law $S$. By
Lemma~\ref{lem:Ct}(1), the law is commutative, has linear part
$(X+Y,U+Z)$, and satisfies $(X,U)*(0,0)=(X,U)$. Associativity of the second
coordinate is Lemma~\ref{lem:Ct}(3). For~(1),
note that $tS(X,Y)=tX+tY+(tX)(tY)$, and that
\[
  U+Z+\bar C_t(X,Y)-\bar H_t(S(X,Y))=\bigl(U-\bar H_t(X)\bigr)
  +\bigl(Z-\bar H_t(Y)\bigr).
\]
The linear part of $\psi$ is $(X,U)\mapsto(tX,\,U-t^{1-p}X)$, which is
invertible over $\Fp[t^{\pm1}]$. For~(2), the addition law of $W_2$ is
\[
  (x_0,x_1)+(y_0,y_1)=\Bigl(x_0+y_0,\ x_1+y_1
  -\frac{(x_0+y_0)^p-x_0^p-y_0^p}{p}\Bigr),
\]
reduced modulo $p$. By Lemma~\ref{lem:Ct}(2) this agrees with the law at
$t=0$: for odd $p$ the signs agree, and for $p=2$ they agree modulo $2$.
For~(3), $L(U)=U$ when $p=2$, so $C_t=t^{-2}(tS-tX-tY)=XY$.
\end{proof}

\subsection{Families defined by structure equations}

\begin{proposition}\label{prop:Gm-family}
Let $r\ge1$. By Theorem~\ref{thm:cartier}(2), there is a smooth formal group
$G$ over $\Fp[[t]]$ whose Cartier module has a $V$-basis
$u,e_1,\dots,e_r$ with
\[
  Fu=[t]u,\qquad Fe_1=u,\qquad Fe_i=e_{i-1}\quad(2\le i\le r).
\]
Its special fibre is isomorphic to $\Wh_{r+1}$, and its geometric generic
fibre is isomorphic to $\Gm\times\Wh_r$.
\end{proposition}

\begin{proof}
By Theorem~\ref{thm:cartier}(3), the fibres have the same $V$-basis and the
specialized structure equations. At $t=0$, put $f_k=e_{r+1-k}$ for
$1\le k\le r$ and $f_{r+1}=u$. Then $Ff_k=f_{k+1}$ for $k\le r$ and
$Ff_{r+1}=0$. These are the structure equations of $\Wh_{r+1}$
(Example~\ref{ex:standard}(b)), so the special fibre is isomorphic to
$\Wh_{r+1}$ by Theorem~\ref{thm:cartier}(2).

Over the algebraically closed field $\kappa(\bar\eta)$, choose $c$ with
$tc^{p-1}=1$, put $a_0=1$, and choose $a_1,\dots,a_r$ recursively with
$ta_i^p=a_{i-1}$. Put $u'=[c]u$ and $e_i'=e_i-[a_i]u$. Their images
$c\bar u$ and $\bar e_i-a_i\bar u$ form a basis of $M/VM$, so they form a
$V$-basis. Using $F[a]=[a^p]F$, we get
\[
  Fu'=[c^pt]u=u',\qquad Fe_1'=u-[ta_1^p]u=0,\qquad
  Fe_i'=e_{i-1}-[ta_i^p]u=e_{i-1}'\quad(2\le i\le r).
\]
These are the structure equations of $\Gm\times\Wh_r$
(Example~\ref{ex:standard}(b), (c)), so the claim follows from
Theorem~\ref{thm:cartier}(2).
\end{proof}

\begin{proposition}\label{prop:unipotent-family}
Let $a\ge b+2$ and $b\ge0$ be integers. By Theorem~\ref{thm:cartier}(2),
there is a smooth formal group $G$ over $\Fp[[t]]$ whose Cartier module has a
$V$-basis $x_1,\dots,x_{a-1},y_1,\dots,y_{b+1}$ with
\[
  Fx_i=x_{i+1}\ (i\le a-2),\quad Fx_{a-1}=0,\quad
  Fy_j=y_{j+1}\ (j\le b),\quad Fy_{b+1}=[t]x_{b+1}.
\]
Its special fibre is isomorphic to $\Wh_{a-1}\times\Wh_{b+1}$, and its
geometric generic fibre is isomorphic to $\Wh_a\times\Wh_b$.
\end{proposition}

\begin{proof}
Note that $x_{b+1}$ exists because $b+1\le a-1$. At $t=0$ the structure
equations consist of two chains, of lengths $a-1$ and $b+1$. As in the
previous proof, the special fibre is $\Wh_{a-1}\times\Wh_{b+1}$.

Over $\kappa(\bar\eta)$, put $g_i=F^{i-1}y_1$ for $1\le i\le a$. Then
$g_i=y_i$ for $i\le b+1$, and $g_i=[t^{p^{i-b-2}}]x_{i-1}$ for
$b+2\le i\le a$. Moreover $Fg_a=0$. If $b\ge1$, choose $s$ with $s^{p^b}=t^{-1}$
and put $z=x_1-[s]y_2$. For $0\le k\le b-1$ we have
$F^kz=x_{k+1}-[s^{p^k}]y_{k+2}$. Since $F^by_2=Fy_{b+1}=[t]x_{b+1}$,
\[
  F^bz=x_{b+1}-[s^{p^b}t]x_{b+1}=0 .
\]
Put $z_k=F^{k-1}z$ for $1\le k\le b$. In $M/VM$, the images of
$g_1,\dots,g_a$ span the subspace $\Sigma$ with basis
$\bar y_1,\dots,\bar y_{b+1},\bar x_{b+1},\dots,\bar x_{a-1}$, and
$\bar z_k\equiv\bar x_k$ modulo $\Sigma$. Hence $g_1,\dots,g_a,z_1,\dots,z_b$
is a $V$-basis. Its structure equations, $Fg_i=g_{i+1}$ for $i<a$,
$Fg_a=0$, $Fz_k=z_{k+1}$ for $k<b$ and $Fz_b=0$, are those of
$\Wh_a\times\Wh_b$. For $b=0$ the elements $g_1,\dots,g_a$ already form a
$V$-basis.
\end{proof}

\begin{proof}[Proof of Theorem~\ref{thm:C}]
Part~(1) follows from Propositions~\ref{prop:Gm-family}
and~\ref{prop:polynomial-law}; for the latter, extend scalars from $\Fp[t]$
to $\Fp[[t]]$. Part~(2) is Proposition~\ref{prop:unipotent-family}.
\end{proof}

\begin{remark}\label{rem:consistency}
These families are consistent with Theorem~\ref{thm:A}. In
Proposition~\ref{prop:Gm-family}, $\dim G[p^m]$ is $\min(m,r)$ generically
and $\min(m,r+1)$ on the special fibre, while the extended polygons are
$\max(y-1,0)$ and $y$. In Proposition~\ref{prop:unipotent-family}, the
dimensions $\min(m,a)+\min(m,b)$ and $\min(m,a-1)+\min(m,b+1)$ agree unless
$b<m<a$, in which case the second is larger by one.
\end{remark}

\begin{remark}\label{rem:dominance}
Let $G_1$ and $G_2$ be unipotent formal groups of dimension $D$ over an
algebraically closed field, with partitions of Witt lengths
$\rho_1=\rho(G_1)$ and $\rho_2=\rho(G_2)$. Then
$\dim G_\nu[p^m]=\sum_i\min(m,\rho_{\nu,i})$ is the sum of the first $m$
parts of the conjugate partition $\rho_\nu'$. Hence the hypotheses of
Question~\ref{q:converse} say that $\rho'_1$ is dominated by $\rho'_2$.
Since conjugation reverses the dominance order \cite{Bry73}, this means that
$\rho_2$ is dominated by $\rho_1$. By \cite[Prop.~2.3]{Bry73}, if $\rho_1$
covers $\rho_2$ in the dominance order, then $\rho_2$ is obtained from
$\rho_1$ by replacing two parts $a$ and $b$, with $a\ge b+2$ and $b\ge0$, by
$a-1$ and $b+1$. Hence the product of the family of
Proposition~\ref{prop:unipotent-family} with a constant unipotent group
realizes every covering relation. An affirmative answer to
Question~\ref{q:converse} for unipotent groups would therefore follow if
specializations could be composed.
\end{remark}

\appendix

\section{\texorpdfstring{The family $X+Y+tXY$}{The family X+Y+tXY}}\label{app:lattice}

Let $G_t$ be the smooth formal group over $R$ with law $X+Y+tXY$. Its special
fibre is $\Ga$, and over $K$ the map $X\mapsto tX$ is an isomorphism from
$G_{t,K}$ onto $\Gm$ with law $x+y+xy$. Let
$J=\ker\bigl(W(R)\to W(k)\bigr)$ and
\[
  I=\ker\bigl(W(R)\to W(R/tR)\bigr),
\]
the Witt vectors over $R$ all of whose components lie in $tR$. For
$w\in W(K)$ with components $(w_0,w_1,\dots)$, let $\gamma_w$ be the curve
\[
  \gamma_w(T)=t^{-1}\Bigl(\prod_{i\ge0}E\bigl(w_iT^{p^i}\bigr)-1\Bigr)
\]
of $G_{t,K}$.

\begin{proposition}\label{prop:lattice}
\begin{enumerate}[label=(\arabic*)]
\item The map $w\mapsto\gamma_w$ is an isomorphism $W(K)\cong\TC(G_{t,K})$
  of $W(K)$-modules. It maps $I$ onto $\TC(G_t)\subseteq\TC(G_{t,K})$, so
  $\TC(G_t)\cong I$ as $W(R)$-modules.
\item The natural map $I\otimes_{W(R)}W(K)\to\TC(G_{t,K})$ is an
  isomorphism.
\item The natural map $I\otimes_{W(R)}W(k)\to\TC(G_{t,k})=\TC(\Ga)$ is
  surjective but not injective.
\end{enumerate}
\end{proposition}

\begin{proof}
(1) The first assertion follows from \eqref{eq:AH} over $K$ and the
isomorphism $X\mapsto tX$. A $p$-typical curve of $G_{t,K}$ comes from
$\TC(G_t)$ if and only if its coefficients lie in $R$. If $w_i\in tR$ for all
$i$, then $\prod_iE(w_iT^{p^i})\in1+tR[[T]]$, since $E$ has coefficients in
$\Zpl$ and constant term $1$. Hence $\gamma_w$ is integral. Conversely, only
the factors with $j\le i$ contribute to the coefficient of $T^{p^i}$ in
$\prod_jE(w_jT^{p^j})$. The factor $j=i$ contributes only through its linear
term, so this coefficient is $w_i+\Pi_i(w_0,\dots,w_{i-1})$ for a polynomial
$\Pi_i$ without constant term. If $\gamma_w$ is integral, induction on $i$
gives $w_i\in tR$ for all $i$, that is, $w\in I$.

(2) Since $[t]\in I$, we have $I[1/[t]]=W(R)[1/[t]]$. As $[t]$ is invertible
in $W(K)$, it follows that $I\otimes_{W(R)}W(K)\cong W(K)$. By (1), this is
the natural map.

(3) The source is $I/JI$. Let $w\in I$. All terms of
$\prod_iE(w_iT^{p^i})-1$ other than $\sum_iw_iT^{p^i}$ are products of at
least two of the $w_j$, hence lie in $t^2R[[T]]$. So the image of the class of
$w$ is $\sum_i(w_i/t \bmod \fm)\,T^{p^i}$. Every curve of $\Ga$ over $k$
arises in this way, by choosing $w_i\in tR$; this proves surjectivity.

Now let $w\in I$ have all components equal to $t^2$. Its image is zero. We
show $w\notin JI$. Choose $q\ge1$ with $p^q>2$.
\begin{enumerate}[label=(\roman*)]
\item Every element of $JI$ lies in $[a]W(R)+p^qI$ for some nonzero
  $a\in\fm$. Indeed, write it as a finite sum $\sum_\nu z_\nu y_\nu$ with
  $z_\nu\in J$ and $y_\nu\in I$, and write
  $z_\nu=\sum_{i<q}p^i[z_{\nu,i}]+p^qz_\nu'$ with $z_{\nu,i}\in\fm$. Since
  $R$ is a valuation ring, the finitely many elements $z_{\nu,i}$ have a
  common divisor $a\in\fm\setminus\{0\}$.
\item The image of $p^qI$ in $W(R/t^{p^q}R)$ is zero. Indeed, for $y\in I$
  the components of $p^qy=V^q\sigma^q(y)$ are $0$ and the $p^q$-th powers
  of the components of $y$, which lie in $t^{p^q}R$.
\item Every element of the image of $[a]W(R)$ in $W(R/t^{p^q}R)$ has only
  finitely many nonzero components. Indeed,
  $[a](x_0,x_1,\dots)=(ax_0,a^px_1,a^{p^2}x_2,\dots)$, and
  $v(a^{p^i})=p^iv(a)\ge p^q$ for large $i$.
\end{enumerate}
The image of $w$ in $W(R/t^{p^q}R)$ has all components equal to $t^2\neq0$.
By (i)--(iii), $w\notin JI$.
\end{proof}

\begin{remark}\label{rem:lattice-generators}
The ideal $I$ is the $p$-adic closure of the ideal $\mathfrak a$ generated by
the elements ${}^{V^i}[t]=p^i[t^{1/p^i}]$, $i\ge0$, but $I\neq\mathfrak a$.
Indeed, an element of $I$ can be written $\sum_ip^i[w_i^{1/p^i}]$ with
$w_i\in tR$, and $p^i[w_i^{1/p^i}]\in p^i[t^{1/p^i}]W(R)$; the truncations of
this sum lie in $\mathfrak a$. On the other hand, every element of
$\mathfrak a$ lies in the ideal generated by $p^i[t^{1/p^i}]$ for $i\le N$,
for some $N$, and this ideal is contained in $[t^{1/p^N}]W(R)$. By the
uniqueness of Teichm\"uller expansions, the element
$\sum_{i\ge0}p^i[t^{1/p^i}]\in I$ does not lie in $[t^{1/p^N}]W(R)$ for any
$N$.
\end{remark}

\section*{Declarations}

\paragraph{Funding.}
Y.Y.~is partially supported by a grant from the National Natural Science
Foundation of China (NSFC-12231001).

\paragraph{Conflict of interest.}
The author declares that there is no conflict of interest.

\paragraph{Data availability.}
Data sharing is not applicable to this article, as no datasets were generated
or analysed.

\paragraph{Use of AI tools.}
The research notes from which this paper developed were drafted with the
assistance of \NotesAI. Claude (Anthropic), with the Claude Opus 5.5 model,
was then used to check the arguments and the references, to reorganize the
notes into the present self-contained article, and to revise the exposition
and the written form of the proofs. It was also used to verify the
polynomial identities of Lemma~\ref{lem:Ct} by computer for $p\le7$. The
author has checked all arguments and takes full responsibility for the
content of this paper.


\begin{thebibliography}{BMS18}

\bibitem[BMS18]{BMS18}
B. Bhatt, M. Morrow and P. Scholze,
\emph{Integral $p$-adic Hodge theory},
Publ. Math. IH\'ES \textbf{128} (2018), 219--397.
\href{https://doi.org/10.1007/s10240-019-00102-z}{doi:10.1007/s10240-019-00102-z}.

\bibitem[Bry73]{Bry73}
T. Brylawski,
\emph{The lattice of integer partitions},
Discrete Math. \textbf{6} (1973), 201--219.
\href{https://doi.org/10.1016/0012-365X(73)90094-0}{doi:10.1016/0012-365X(73)90094-0}.

\bibitem[Car67]{Car67}
P. Cartier,
\emph{Modules associ\'es \`a un groupe formel commutatif. Courbes typiques},
C. R. Acad. Sci. Paris S\'er. A-B \textbf{265} (1967), A129--A132.

\bibitem[Dem72]{Dem72}
M. Demazure,
\emph{Lectures on $p$-divisible groups},
Lecture Notes in Mathematics, vol.~302, Springer, Berlin, 1972.
\href{https://doi.org/10.1007/BFb0060741}{doi:10.1007/BFb0060741}.

\bibitem[Gro74]{Gro74}
A. Grothendieck,
\emph{Groupes de Barsotti--Tate et cristaux de Dieudonn\'e},
S\'eminaire de Math\'ematiques Sup\'erieures, no.~45 (\'Et\'e 1970),
Les Presses de l'Universit\'e de Montr\'eal, Montreal, 1974.

\bibitem[Haz78]{Haz78}
M. Hazewinkel,
\emph{Formal groups and applications},
Pure and Applied Mathematics, vol.~78, Academic Press, New York--London,
1978.

\bibitem[Kat79]{Kat79}
N. M. Katz,
\emph{Slope filtration of $F$-crystals},
in: Journ\'ees de G\'eom\'etrie Alg\'ebrique de Rennes (Rennes, 1978),
Vol.~I, Ast\'erisque \textbf{63} (1979), 113--163.

\bibitem[Man63]{Man63}
Yu. I. Manin,
\emph{The theory of commutative formal groups over fields of finite
characteristic},
Russian Math. Surveys \textbf{18} (1963), no.~6, 1--83.

\bibitem[Oor00]{Oor00}
F. Oort,
\emph{Newton polygons and formal groups: conjectures by Manin and
Grothen\-dieck},
Ann. of Math. (2) \textbf{152} (2000), 183--206.

\bibitem[Zin84]{Zin84}
T. Zink,
\emph{Cartiertheorie kommutativer formaler Gruppen},
Teubner-Texte zur Mathematik, vol.~68, Teubner, Leipzig, 1984.
English translation by M. Garuti, M. Le Barbier Gr\"unewald, C. P\'epin
and M. Romagny: \emph{Cartier theory of commutative formal groups}.

\end{thebibliography}
\end{document}